\documentclass[pdflatex,sn-mathphys-num]{sn-jnl}% Math and Physical Sciences Numbered Reference Style
\usepackage{graphicx}%
\usepackage{multirow}%
\usepackage{amsmath,amssymb,amsfonts,mathtools}%
\usepackage{amsthm}%
\usepackage{mathrsfs}%
\usepackage[title]{appendix}%
\usepackage{xcolor}%
\usepackage{textcomp}%
\usepackage{manyfoot}%
\usepackage{booktabs}%
\usepackage{algorithm}%
\usepackage{algorithmicx}%
\usepackage{algpseudocode}%
\usepackage{listings}%
\usepackage{bm}
\usepackage{microtype}
\usepackage{enumitem}
\usepackage{dsfont}

\theoremstyle{thmstyleone}%
\newtheorem{theorem}{Theorem}%  meant for continuous numbers
\newtheorem{proposition}[theorem]{Proposition}% 

\theoremstyle{thmstyletwo}%
\newtheorem{corollary}[theorem]{Corollary}%
\DeclareMathOperator{\Var}{Var}
\newcommand{\R}{\mathbb{R}}
\newcommand{\rev}{\color{black}}
\newcommand{\revv}{\color{black}} % use \relax instead of \color{red} to print revisions in black

\theoremstyle{thmstylethree}%
\begin{document}

\title[Decision-Weighted Limit Theorems: Less Conservative Robust Optimization]{Decision-Weighted Limit Theorems: Less Conservative Robust Optimization}

%%=============================================================%%
%% GivenName	-> \fnm{Joergen W.}
%% Particle	-> \spfx{van der} -> surname prefix
%% FamilyName	-> \sur{Ploeg}
%% Suffix	-> \sfx{IV}
%% \author*[1,2]{\fnm{Joergen W.} \spfx{van der} \sur{Ploeg} 
%%  \sfx{IV}}\email{iauthor@gmail.com}
%%=============================================================%%

\author[1]{\fnm{Dimitris} \sur{Bertsimas}}\email{dbertsim@mit.edu}

\author*[1]{\fnm{Danique} \sur{de Moor}}\email{daniquedemoor@hotmail.com}
% \equalcont{These authors contributed equally to this work.}

\author[2]{\fnm{Dick} \sur{den Hertog}}\email{d.denhertog@uva.nl}
% \equalcont{These authors contributed equally to this work.}

\author[3]{\fnm{Federica} \sur{Donnini}}\email{federica.donnini@unifi.it}
% \equalcont{These authors contributed equally to this work.}

\affil[1]{\orgdiv{Sloan School of Management}, \orgname{MIT}, \city{Cambridge},  \country{USA}}

\affil[2]{\orgdiv{Amsterdam Business School}, \orgname{University of Amsterdam},  \city{Amsterdam}, \country{The Netherlands}}

\affil[3]{\orgdiv{Department of Information Engineering}, \orgname{University of Florence}, \orgaddress{\city{Florence}, \country{Italy}}}

%%==================================%%
%% Sample for unstructured abstract %%
%%==================================%%
\abstract{Limit theorems provide a natural way to construct tractable models for optimization under uncertainty. Existing approaches typically apply them to unweighted sums of uncertain parameters and use the resulting bounds to construct uncertainty sets. We show that sharper formulations are obtained by applying limit theorems directly to the uncertain linear expressions in the optimization problem, with the decision variables as weights. For independent random variables with finite variance, the weighted central limit theorem yields a second-order cone constraint with a decision-dependent protection term, and for heavy-tailed random variables governed by stable limit laws, the weighted approach yields power-cone constraints. Using norm duality, we show that, under a condition on the individual bounds, the feasible region of the unweighted formulation is contained in that of the weighted one. In 16,000 knapsack instances with a prescribed violation probability of 5\%, the weighted formulation improves the optimal objective value by 20.9\% on average, with overall average violation probabilities between 3.45\% and 4.95\%, whereas those of the unweighted formulation are near zero.}


%%\pacs[JEL Classification]{D8, H51}


\maketitle

\section{Introduction}
Limit theorems provide a natural bridge between probability theory and optimization under uncertainty. In particular, Bandi and Bertsimas \cite{bandi-bertsimas} propose using asymptotic laws of probability to construct uncertainty sets for robust optimization. For example, the Central Limit Theorem (CLT) provides bounds on aggregate sums of independent random variables, which can then be incorporated into uncertainty sets. This approach has been used in inventory management \cite{mamani,sun-vanmieghem,thorsen-yao,xie,Gijsbrechts2023}, queueing systems \cite{BandiBertsimasYoussef2015}, option pricing \cite{BandiBertsimas2014}, and transportation procurement \cite{Zhang2015}. CLT-based uncertainty sets are also discussed in the textbook by Bertsimas and den Hertog \cite{BertsimasDenHertog2022}. These sets require limited distributional information and often yield tractable robust optimization models. However, their construction typically applies limit theorems to unweighted aggregates, without exploiting decision-dependent weights. While natural when the aggregate itself is of interest, this may introduce conservatism for decision-weighted sums. This motivates the weighted limit-theorem approach developed in this paper.

% Limit theorems provide a natural bridge between probability theory and optimization under uncertainty. In particular, Bandi and Bertsimas~\cite{bandi-bertsimas} propose using asymptotic laws of probability to construct uncertainty sets for robust optimization. For example, the Central Limit Theorem (CLT) provides bounds on the aggregate sum of a collection of independent random variables, which can then be incorporated into an uncertainty set. Such sets require limited distributional information, lead to tractable models, and have been adopted in several applications ~\cite{mamani, sun-vanmieghem, thorsen-yao, xie} ({\color{red} cite also other references here not in robust inventory management}).

In optimization, uncertain coefficients typically enter through expressions of the form $\bm a^\top\bm x$, where $\bm a\in\R^n$ collects the uncertain coefficients and $\bm x\in\R^n$ the decision variables. In existing constructions~\cite{bandi-bertsimas}, the CLT is applied to the aggregate sum $\sum_j a_j$, and the resulting bound defines an uncertainty set that protects $\bm a^\top\bm x$. This ignores the decision weights $x_j$: the aggregate bound is independent of $\bm x$, although the distribution of $\bm a^\top\bm x$ may depend strongly on $\bm x$. The unweighted construction nevertheless has advantages: Bandi and Bertsimas~\cite{bandi-bertsimas} consider settings, such as queueing networks, where the aggregate itself is of interest, and their decision-independent uncertainty sets can be reused across constraints.

In this note, we apply limit theorems directly to the weighted sum $\bm a^\top\bm x$, treating the decision variables $x_j$ as weights. The resulting constraints reflect the scaling behavior of the underlying limit theorem: the weighted CLT yields second-order cone constraints, while weighted stable limit theorems yield power-cone constraints.

Second-order cone reformulations of chance constraints with jointly Gaussian coefficients are classical~\cite{nemirovski-shapiro}. Our approach does not require normally distributed coefficients: normality arises asymptotically from aggregating the independent terms $x_j a_j$, yielding the same conic structure for a broader class of distributions. The guarantee is therefore asymptotic rather than finite-sample. In contrast, safe approximations, such as the Bernstein approximation~\cite{nemirovski-shapiro} and moment-based distributionally robust bounds~\cite{elghaoui2003worstcase,zymler2013distributionally}, provide finite-sample guarantees under their respective assumptions, but are generally more conservative. Moreover, moment-based bounds do not fully exploit independence.

Because the guarantee is asymptotic and the weights are selected by optimization, we bound the normal-approximation error at each decision in terms of the largest term's contribution relative to the standard deviation of $\bm a^\top\bm x$, yielding a reliability guarantee for the selected decision. Using norm duality, we further show that, under a condition on the individual bounds, every decision feasible for the unweighted formulation is also feasible for the weighted one. In knapsack experiments, the weighted formulation improves the optimal objective value by $20.9\%$ on average for a prescribed maximum violation probability of $\varepsilon=0.05$, with average violation probabilities of $3.45\%$--$4.95\%$ across three distributions, compared with near-zero probabilities for the unweighted formulation.

\section{Weighted Central Limit Theorem}\label{sec:weighted CLT}
Consider an optimization problem containing an uncertain linear constraint $\bm{a}^\top\bm{x}\le b$, where $\bm{x}\in\mathcal X_n\subseteq\R^n$ and $\bm{a}\in\R^n$ consists of independent random variables with $\mathbb E[a_j]=\mu_j$ and $\Var(a_j)=\sigma_j^2>0$. {\rev Apart from independence and finite variances,} no distributional assumption is imposed on the coefficients $a_j$. {\rev An uncertain linear objective can be treated in the same way using the standard epigraph reformulation. }
% Since limit theorems describe the regime $n\to\infty$, we consider a sequence of such problems indexed by $n$, in which $a_j$, $\mu_j$, $\sigma_j$, $b$, and $\mathcal X_n$ may depend on $n$; this dependence is suppressed in the notation.
% Since limit theorems describe the regime $n\to\infty$, we consider a sequence of such problems indexed by $n$, in which $a_j$, $\mu_j$, $\sigma_j$, $b$, and $\mathcal X_n$ may depend on $n$. This dependence is suppressed in the notation.

We aim to approximate the chance constraint
\begin{equation}\label{eq:chance constraint}
    \mathbb P(\bm{a}^\top\bm{x}\le b)\ge 1-\varepsilon,
\end{equation}
where $\varepsilon\in(0,1/2)$ is the prescribed violation probability. Throughout, $\Phi$ denotes the cumulative distribution function of the standard normal distribution $\mathcal N(0,1)$, $z_{1-\varepsilon}=\Phi^{-1}(1-\varepsilon)$, $\bm D=\operatorname{diag}(\sigma_1,\ldots,\sigma_n)$, and $\xrightarrow{d}$ denotes convergence in distribution. %~as~$n\to\infty$.

%\medskip
%\noindent \textbf{Unweighted construction.}
\bmhead{Unweighted construction}
Following Bandi and Bertsimas~\cite{bandi-bertsimas}, the CLT can be applied to the unweighted sum of the uncertain coefficients. Under standard conditions,
\[
    \frac{\sum_{j=1}^n(a_j-\mu_j)}{\sqrt{\sum_{j=1}^n\sigma_j^2}}\xrightarrow{d}\mathcal N(0,1),
\]
which motivates an aggregate bound of the form
\[
    \Bigl|\sum_{j=1}^n(a_j-\mu_j)\Bigr|\le  C\Bigg(\sum_{j=1}^n\sigma_j^2\Bigg)^{1/2},\qquad C>0.
\]
{\rev Since this bound is two-sided, it holds with asymptotic probability $1-\varepsilon$ for $C=z_{1-\varepsilon/2}$.
Combining it with individual bounds
$|a_j-\mu_j|\leq\eta_j$ for each $j$, yields the
uncertainty set
\[
    \mathcal U_U^{(2)}=\bigl\{\bm a\in\R^n:\ |a_j-\mu_j|\le\eta_j,\ j=1,\ldots,n,\ \ |\bm 1^\top(\bm a-\bm\mu)|\le C\lVert\bm\sigma\rVert_2\bigr\},
\]
where $\bm 1$ is the all-ones vector.}

\bmhead{Weighted construction}
%\medskip
%\noindent \textbf{Weighted construction.}
We instead apply the CLT directly to the decision-weighted expression $\bm{a}^\top\bm{x}$. For $\bm{x}\in\mathcal X_n$, define $S(\bm{x})=\sum_{j=1}^n x_j(a_j-\mu_j)$. By independence, $\mathbb E[S(\bm{x})]=0$ and $\Var(S(\bm{x}))=\lVert\bm D\bm{x}\rVert_2^2$. {\revv 
For a sequence of nonzero decision vectors
$\bm x^{(n)}\in\mathcal X_n$ as $n\to\infty$,
if the corresponding weighted terms satisfy
the Lindeberg condition, the CLT gives}
\begin{equation}\label{eq:conditions weighted CLT}
    \frac{S(\bm{x})}{\lVert\bm D\bm{x}\rVert_2}\xrightarrow{d}\mathcal N(0,1).
\end{equation}
{\rev Since $\mathbb P(\bm{a}^\top\bm{x}\le b)=\mathbb P\bigl(S(\bm{x})/\lVert\bm D\bm{x}\rVert_2\le(b-\bm{\mu}^\top\bm{x})/\lVert\bm D\bm{x}\rVert_2\bigr)$, replacing the distribution of the standardized sum by its normal limit in~\eqref{eq:chance constraint} yields the asymptotic deterministic constraint
\begin{equation}\label{eq:asymptotic deterministic constraint equiv}
    \bm{\mu}^\top\bm{x}+z_{1-\varepsilon}\lVert \bm D\bm{x}\rVert_2\le b.
\end{equation}
Since $\varepsilon<1/2$, we have $z_{1-\varepsilon}>0$, so~\eqref{eq:asymptotic deterministic constraint equiv} is a second-order cone constraint.}

The key difference with the unweighted construction is that the protection term depends on the decision: the relevant standard deviation is that of $\bm{a}^\top\bm{x}$, namely $\lVert\bm D\bm{x}\rVert_2=(\sum_{j=1}^n\sigma_j^2x_j^2)^{1/2}$, rather than that of the unweighted sum, namely $(\sum_{j=1}^n\sigma_j^2)^{1/2}$. The weighted construction therefore uses the uncertainty of precisely the linear combination appearing in the constraint, whereas an uncertainty set derived from the unweighted sum is calibrated without knowing which linear combination will be relevant. {\rev We investigate the implications of this distinction in Sections~\ref{sec:theoretical comparison} and~\ref{sec:numerical comparison}.}
\\
\\
\noindent \textbf{Uniformity over decisions.}
{\revv The normal approximation in~\eqref{eq:asymptotic deterministic constraint equiv} relies on aggregation, not merely on a large dimension $n$. When the variance of $\bm a^\top\bm x$ is spread across many independent terms, the sum approaches normality, whereas when concentrated in a few, deviations from normality may persist regardless of $n$. 
In the extreme case $\bm x=\bm e_j$, where $\bm e_j\in\mathbb R^n$
denotes the $j$-th unit vector, the approximation treats
a single coefficient $a_j$ as normally distributed.

The weighted CLT~\eqref{eq:conditions weighted CLT} guarantees that the approximation error vanishes along sequences of sufficiently spread-out decisions, but provides no error bound at a given decision. Such a bound is important in optimization, where the selected decision may concentrate its weight on only a few terms. We measure this error at $\bm x\neq\bm 0$ by
\[
    \delta(\bm x)=\sup_{t\in\R}\left|\mathbb P\left(\frac{S(\bm{x})}{\lVert\bm D\bm{x}\rVert_2}\le t\right)-\Phi(t)\right|,
\]
the largest discrepancy between the true and approximated violation probabilities over all thresholds $b$. We bound it in terms of
\begin{equation}\label{eq:cond_nodominantcoefficients}
    R(\bm{x})=\frac{\max_{j=1,\ldots,n}\sigma_j|x_j|}{\lVert\bm D\bm{x}\rVert_2},
    \qquad \bm{x}\in\R^n\setminus\{\bm{0}\},
\end{equation}
which measures the largest individual contribution to the standard deviation of the sum. 
This ratio lies in $[n^{-1/2},1]$, equals $1$ when a single term accounts for all the variance, and can be evaluated at the selected decision. 
For $r\in(0,1]$, define
\begin{equation}\label{eq:BE-bound}
    \psi(r)=\sup_{n}\ \sup_{\bm{x}\in\R^n\setminus\{\bm 0\}:\,R(\bm x)\le r}\ \delta(\bm x),
\end{equation}
the worst-case approximation error over all dimensions and decisions with $R(\bm x)\le r$. The following proposition shows that $\psi(r)\to0$ as $r\to0$, uniformly in $n$.

\begin{proposition}[Uniform weighted CLT]\label{prop:uniform-clt}
Let $Y_j=(a_j-\mu_j)/\sigma_j$ and suppose that the $Y_j^2$ are uniformly integrable, i.e.,
\begin{equation}\label{eq:lyapunov}
    \lim_{K\to\infty}\ \sup_{n}\ \max_{j=1,\ldots,n}\ \mathbb E\bigl[Y_j^2\,\mathds{1}_{\{|Y_j|>K\}}\bigr]=0,
\end{equation}
where $\mathds{1}_A$ denotes the indicator of the event $A$. Then $\psi(r)\to0$ as $r\to0$.
\end{proposition}

\begin{proof}
For $\bm x\neq\bm 0$, let $w_j=\sigma_jx_j/\lVert\bm D\bm x\rVert_2$. Then $S(\bm x)/\lVert\bm D\bm x\rVert_2=\sum_{j=1}^{n}w_jY_j$, $\sum_{j=1}^{n}w_j^2=1$, and $\max_{j=1,\ldots,n}|w_j|=R(\bm x)$.
Suppose that $\psi(r)\not\to0$. Since $\psi$ is nondecreasing, there exists $\eta>0$ such that $\psi(r)>\eta$ for all $r\in(0,1]$. By the definition of $\psi$ with $r=1/k$, there exist, for every $k\ge1$, a problem size $n_k$ and a decision $\bm x^{(k)}\in\R^{n_k}\setminus\{\bm 0\}$ with $R(\bm x^{(k)})\le1/k$ and $\delta(\bm x^{(k)})>\eta$. Let $w_{k,j}$ and $Y_{k,j}$ denote $w_j$ and $Y_j$ for the decision $\bm x^{(k)}$ and the problem of size $n_k$. Then $\sum_{j=1}^{n_k}w_{k,j}^2=1$, $\max_j|w_{k,j}|\le1/k$, and
\begin{equation}\label{eq:bad-sequence}
    \sup_{t}\Bigl|\mathbb P\Bigl(\sum_{j}^{n_k}w_{k,j}Y_{k,j}\le t\Bigr)-\Phi(t)\Bigr|>\eta\qquad\text{for all }k.
\end{equation}
Let $X_{k,j}=w_{k,j}Y_{k,j}$. Since the $Y_{k,j}$ are independent with mean zero, unit variance, and $\sum_{j=1}^{n_k}w_{k,j}^2=1$, the variables $X_{k,1},\ldots,X_{k,n_k}$ are independent with mean zero, and their variances sum to one. We show that the triangular array with rows $X_{k,1},\ldots,X_{k,n_k}$, $k=1,2,\ldots$, satisfies the Lindeberg condition, that is, for every $\tau>0$,
\[
    L_k(\tau)=\sum_{j=1}^{n_k}\mathbb E\bigl[X_{k,j}^2\,\mathds 1_{\{|X_{k,j}|>\tau\}}\bigr]\longrightarrow0\qquad(k\to\infty).
\]
Since $|w_{k,j}|\le1/k$, the event $\{|X_{k,j}|>\tau\}$ is contained in $\{|Y_{k,j}|>k\tau\}$. With $g(K)=\sup_n\max_{j=1,\ldots,n}\mathbb E[Y_j^2\,\mathds 1_{\{|Y_j|>K\}}]$, it follows that
\[  L_k(\tau)\le\sum_{j=1}^{n_k}w_{k,j}^2\,\mathbb E\bigl[Y_{k,j}^2\,\mathds 1_{\{|Y_{k,j}|>k\tau\}}\bigr]\le g(k\tau)\sum_{j=1}^{n_k}w_{k,j}^2=g(k\tau),
\]
which tends to zero as $k\to\infty$. This bound depends on the weights only through their largest absolute value, which is why it holds uniformly over all decisions.

By the Lindeberg--Feller CLT (see, e.g.,~\cite{petrov1995}), $\sum_{j=1}^{n_k}X_{k,j}\xrightarrow{d}\mathcal N(0,1)$ as $k\to\infty$, that is, $\mathbb P(\sum_{j=1}^{n_k}X_{k,j}\le t)\to\Phi(t)$ for every $t$. Since $\Phi$ is continuous, P\'olya's theorem (see, e.g.,~\cite[Lemma~2.11]{vandervaart1998}) implies that this convergence is uniform in $t$, which contradicts~\eqref{eq:bad-sequence}. Hence,~$\psi(r)\to0$.
\end{proof}
Since $\delta(\bm x)$ bounds the error that~\eqref{eq:asymptotic deterministic constraint equiv} makes in the violation probability, Proposition~\ref{prop:uniform-clt} yields the following reliability guarantee, which applies to every decision satisfying~\eqref{eq:asymptotic deterministic constraint equiv}.
\begin{corollary}\label{cor:violation}
Every $\bm{x}\in\mathcal X_n\setminus\{\bm 0\}$ satisfying~\eqref{eq:asymptotic deterministic constraint equiv} satisfies
$\mathbb P(\bm{a}^\top\bm{x}>b)\le\varepsilon+\psi(R(\bm x)).$
% \begin{equation*}
%         \mathbb P(\bm{a}^\top\bm{x}>b)\le\varepsilon+\psi(R(\bm x)).
% \end{equation*}
\end{corollary}
\begin{proof}
From~\eqref{eq:asymptotic deterministic constraint equiv} it holds $t=(b-\bm{\mu}^\top\bm{x})/\lVert\bm D\bm{x}\rVert_2\ge z_{1-\varepsilon}$. Since $\mathbb P(\bm{a}^\top\bm{x}\le b)=\mathbb P(S(\bm{x})/\lVert\bm D\bm{x}\rVert_2\le t)$, the definition of $\psi$ gives $\mathbb P(\bm{a}^\top\bm{x}>b)\le1-\Phi(t)+\psi(R(\bm x))\le\varepsilon+\psi(R(\bm x))$.
\end{proof}
\section{Weighted Stable Limit Theorems}\label{sec:weighted stable limit thms}
The CLT of Section~\ref{sec:weighted CLT} requires finite variances. In many applications, however, uncertain quantities are heavy-tailed, and their variances may not exist. As discussed by Bandi and Bertsimas~\cite{bandi-bertsimas}, stable limit laws provide an analogue of the CLT for such distributions and can therefore also be used to construct tractable models for optimization under uncertainty. In this section, we show that the same principle as in Section~\ref{sec:weighted CLT} applies: rather than applying a stable limit theorem to the unweighted sum of the uncertain coefficients, we apply it directly to the decision-weighted expression $\bm{a}^\top\bm{x}$.

We first consider coefficients that are themselves stable. Let $a_j=\mu_j+\gamma_jZ_j$, $j=1,\ldots,n$, {\rev where $\gamma_j>0$ is a scale parameter} and the $Z_j$ are independent standard symmetric $\alpha$-stable random variables, {\rev that is, $\mathbb E[e^{\mathrm itZ_j}]=e^{-|t|^\alpha}$ for $t\in\R$,} with common stability parameter $1<\alpha<2$. The restriction $\alpha>1$ ensures that the mean exists and that $\lVert\cdot\rVert_\alpha$ is a norm, while $\alpha<2$ corresponds to the heavy-tailed case with infinite variance. For background on stable distributions, see, e.g.,~\cite{samorodnitsky-taqqu}.

A defining property of stable distributions is that weighted sums remain stable. Let $\overset{d}{=}$ denote equality in distribution, then for a fixed vector $\bm{x}$, it holds
\begin{equation*}
    \sum_{j=1}^n\gamma_jx_jZ_j
    \overset{d}{=}
    \Big(\sum_{j=1}^n|\gamma_jx_j|^\alpha\Big)^{1/\alpha}Z\quad \Rightarrow\quad \bm{a}^\top\bm{x}
    \overset{d}{=}
    \bm{\mu}^\top\bm{x}
    +\Big(\sum_{j=1}^n|\gamma_jx_j|^\alpha\Big)^{1/\alpha}Z,
\end{equation*}
{\rev where $Z$ is a standard symmetric $\alpha$-stable random variable}. %Consequently,
% \begin{equation*}\label{eq:stable-weighted-sum}
%     \bm{a}^\top\bm{x}
%     \overset{d}{=}
%     \bm{\mu}^\top\bm{x}
%     +\Big(\sum_{j=1}^n|\gamma_jx_j|^\alpha\Big)^{1/\alpha}Z.
% \end{equation*}
%Let $q_{1-\varepsilon}^{(\alpha)}$ {\color{red} be }the $(1-\varepsilon)$-quantile of the distribution {\color{red}of $Z$} {\rev and  $\bm \Gamma=\operatorname{diag}(\gamma_1,\ldots,\gamma_n)$.} Then  chance constraint~\eqref{eq:chance constraint} is equivalent to
{\rev Let $q_{1-\varepsilon}^{(\alpha)}$ denote the $(1-\varepsilon)$-quantile of the standard symmetric $\alpha$-stable distribution and  $\bm \Gamma=\operatorname{diag}(\gamma_1,\ldots,\gamma_n)$.} Then~\eqref{eq:chance constraint}
is equivalent to
% \begin{equation}
%     \bm{\mu}^\top\bm{x}
%     +q_{1-\varepsilon}^{(\alpha)}
%     \left(\sum_{j=1}^n|\gamma_jx_j|^\alpha\right)^{1/\alpha}
%     \le b,
% \end{equation}
% which, in terms of $\Gamma$, reads
\begin{equation}\label{eq:chance constraint eq}
    \bm{\mu}^\top\bm{x}
    +q_{1-\varepsilon}^{(\alpha)}\lVert \bm \Gamma\bm{x}\rVert_\alpha
    \le b.
\end{equation}
{\rev By symmetry of $Z$, $q^{(\alpha)}_{1-\varepsilon}>0$ for $\varepsilon<1/2$, so~\eqref{eq:chance constraint eq} is convex.} Moreover, using an epigraph reformulation, it is power-cone representable. The stable distribution thus leads to a power-cone formulation in the same way that the normal distribution leads to a second-order cone formulation in Section~\ref{sec:weighted CLT}. {\rev The quantile $q^{(\alpha)}_{1-\varepsilon}$ depends only on $\alpha$ and $\varepsilon$. Since the stable distribution function has no closed form for $1<\alpha<2$, it is computed numerically, once per model, using standard software (see, e.g.,~\cite{nolan2020}).}

As in the finite-variance case, the protection term depends on the decision: the unweighted sum $\sum_j(a_j-\mu_j)$ has scale $(\sum_{j}\gamma_j^\alpha)^{1/\alpha}$, independent of $\bm{x}$, whereas $\bm{a}^\top\bm{x}$ has scale $\lVert\Gamma\bm{x}\rVert_\alpha$. The weighted formulation therefore accounts directly for which uncertain coefficients enter the constraint and with what weights.

The assumption that the coefficients are themselves stable can be relaxed. Suppose instead that the centered coefficients $a_j-\mu_j$ belong to the domain of attraction of an $\alpha$-stable law. Under suitable regularity conditions on the distributions and on a sequence of weights, a weighted stable limit theorem yields an asymptotic relation of the form 
\begin{equation*}\label{eq:weighted stable limit}
    \frac{\sum_{j=1}^nx_j(a_j-\mu_j)}{\lVert\bm \Gamma\bm{x}\rVert_\alpha}\xrightarrow{d}Z,
\end{equation*}
so that~\eqref{eq:chance constraint eq} becomes an asymptotic approximation of~\eqref{eq:chance constraint} rather than an exact reformulation. {\rev Among these conditions, no single term may dominate $\bm{a}^\top\bm{x}$, that is,

\[
R_\alpha(\bm{x})=\frac{\max_j\gamma_j|x_j|}{\lVert\bm \Gamma\bm{x}\rVert_\alpha}\longrightarrow0
\]
along the sequence of decisions considered, where $R_\alpha(\bm{x})\in[n^{-1/\alpha},1]$ is the analogue of $R(\bm{x})$ in~\eqref{eq:cond_nodominantcoefficients}. Explicit finite-$n$ error bounds would additionally require rates of convergence in the domain of attraction, which depend on second-order properties of the tails of the $a_j$. When the coefficients are themselves stable,~\eqref{eq:chance constraint eq} is exact and no such condition is needed.}

\section{Comparison with Unweighted Limit-Theorem Uncertainty Sets}\label{sec:theoretical comparison}
We use norm duality to show that, under a
condition on the individual bounds, every decision feasible
for the unweighted formulation is also feasible for the
weighted one, in both the finite-variance and stable settings. Let $p>1$, $C>0$, and $p^*=p/(p-1)$.
The case $p=2$ corresponds to Section~\ref{sec:weighted CLT},
and $p=\alpha$ to Section~\ref{sec:weighted stable limit thms},
with $\sigma_j=\gamma_j$ and $\bm D=\bm\Gamma$.

Generalizing the unweighted uncertainty set
$\mathcal U_U^{(2)}$ from Section 2 to $p>1$, we obtain
\begin{equation*}
\begin{aligned}
\mathcal U_U^{(p)}
=\bigl\{\bm a\in\R^n:\quad
&|a_j-\mu_j|\le\eta_j,\quad j=1,\ldots,n, 
&|\bm 1^\top(\bm a-\bm\mu)|
\le C\lVert\bm\sigma\rVert_p
\bigr\}.
\end{aligned}
\end{equation*}
The weighted construction yields the constraint $\bm\mu^\top\bm x+C\lVert\bm D\bm x\rVert_p\le b$. By the duality between the $p$- and $p^*$-norms and H\"older's inequality, this is the robust counterpart of $\bm a^\top\bm x\le b$ for the uncertainty set
\begin{equation*}
    \mathcal U_W^{(p)}=\bigl\{\bm\mu+\bm D\bm z:\ \lVert\bm z\rVert_{p^*}\le C\bigr\},
\end{equation*}
since
\[
    \max_{\bm a\in\mathcal U_W^{(p)}}\bm a^\top\bm x=\bm\mu^\top\bm x+\max_{\lVert\bm z\rVert_{p^*}\le C}\bm z^\top\bm D\bm x=\bm\mu^\top\bm x+C\lVert\bm D\bm x\rVert_p.
\]
The next result provides a comparison of the weighted and unweighted constructions.
\begin{theorem}\label{thm:dominance}
If $\eta_j\ge C\sigma_j$ for $j=1,\ldots,n$, then $\mathcal U_W^{(p)}\subseteq\mathcal U_U^{(p)}$. Consequently, for every $\bm x\in\R^n$,
\[
    \bm\mu^\top\bm x+C\lVert\bm D\bm x\rVert_p\le\max_{\bm a\in\mathcal U_U^{(p)}}\bm a^\top\bm x.
\]
\end{theorem}
\begin{proof}
%\noindent\emph{Proof}
Let $\bm{a}\in\mathcal U_W^{(p)}$. There exists $\bm{z}\in\mathbb R^n$ such that $\bm{a}=\bm{\mu}+\bm D\bm{z}$ and $\lVert\bm{z}\rVert_{p^*}\le C.$
% \[
%     \bm{a}=\bm{\mu}+\bm D\bm{z},
%     \qquad
%     \lVert\bm{z}\rVert_{p^*}\le C.
% \]
For every $j$,
\[
    |a_j-\mu_j|
    =\sigma_j|z_j|
    \le\sigma_j\lVert\bm{z}\rVert_{p^*}
    \le C\sigma_j
    \le\eta_j.
\]
Hence, all individual bounds defining $\mathcal U_U^{(p)}$ are satisfied. Moreover, by H\"older's inequality,
\[
    |\bm{1}^\top(\bm{a}-\bm{\mu})|
    =|\bm{1}^\top \bm D\bm{z}|
    =|\bm{\sigma}^\top\bm{z}|\le\lVert\bm{\sigma}\rVert_p\lVert\bm{z}\rVert_{p^*}
    \le C\lVert\bm{\sigma}\rVert_p.
\]
%where the first inequality follows by H\"older's inequality.
% \[
%     |\bm{\sigma}^\top\bm{z}|
%     \le\lVert\bm{\sigma}\rVert_p\lVert\bm{z}\rVert_{p^*}
%     \le C\lVert\bm{\sigma}\rVert_p.
% \]
Thus, the aggregate constraint of $\mathcal U_U^{(p)}$ is also satisfied, and
$\bm{a}\in\mathcal U_U^{(p)}$. Hence $\mathcal U_W^{(p)}\subseteq\mathcal U_U^{(p)}$, and
\[
    \max_{\bm{a}\in\mathcal U_W^{(p)}}\bm{a}^\top\bm{x}
    \le
    \max_{\bm{a}\in\mathcal U_U^{(p)}}\bm{a}^\top\bm{x}
\]
follows as a consequence for every $\bm{x}\in\mathbb R^n$. The stated inequality can be obtained using
\[
    \max_{\bm{a}\in\mathcal U_W^{(p)}}\bm{a}^\top\bm{x}
    =\bm{\mu}^\top\bm{x}+C\lVert \bm D\bm{x}\rVert_p.
\]
%\hfill$\square$
\end{proof}
Theorem~\ref{thm:dominance} does not require $\bm x$ to be binary or nonnegative. For a deterministic set $\mathcal X\subseteq\R^n$, define the feasible regions
\[
    \mathcal F_U^{(p)}=\Bigl\{\bm x\in\mathcal X:\max_{\bm a\in\mathcal U_U^{(p)}}\bm a^\top\bm x\le b\Bigr\},\qquad
    \mathcal F_W^{(p)}=\bigl\{\bm x\in\mathcal X:\bm\mu^\top\bm x+C\lVert\bm D\bm x\rVert_p\le b\bigr\}.
\]
\begin{corollary}\label{cor:feasible-regions}
If $\eta_j\ge C\sigma_j$ for $j=1,\ldots,n$, then $\mathcal F_U^{(p)}\subseteq\mathcal F_W^{(p)}$ for every $\mathcal X\subseteq\R^n$ and every $b\in\R$. Consequently, for any objective function, the optimal value of the weighted formulation is at least as good as that of the unweighted formulation.
\end{corollary}
{\rev The condition $\eta_j\ge C\sigma_j$ requires the individual bounds to be no tighter than $C$ standard deviations. Without it, dominance can fail: for $\bm x=\bm e_j$ with $\eta_j<C\sigma_j$, the unweighted protection equals $\eta_j$, whereas the weighted protection equals $C\sigma_j$.}

{\rev The two constructions use different projections of $\bm{a}-\bm{\mu}$. The unweighted construction controls only} the single projection $\bm{1}^\top(\bm{a}-\bm{\mu})$ and combines this information with the individual bounds to protect arbitrary directions $\bm{x}$, whereas the weighted construction applies the limit theorem directly to the projection $\bm{x}^\top(\bm{a}-\bm{\mu})$ that appears in the optimization problem. {\rev Since the aggregate constraint defining $\mathcal U_U^{(p)}$ restricts only the direction $\bm{1}$, $\mathcal U_U^{(p)}$ contains coefficient vectors whose deviations in other directions can far exceed those allowed by $\mathcal U_W^{(p)}$. Specifically, positive deviations of the coefficients with $x_j>0$ can be offset in the aggregate constraint by negative deviations of the remaining coefficients, so that the aggregate constraint does not limit the joint deviation of the coefficients that enter the constraint. 
We refer to this as the \emph{compensation mechanism}. It is the source of the conservatism of the unweighted construction, and it shows that the inclusions in Theorem~\ref{thm:dominance} and Corollary~\ref{cor:feasible-regions} can be strict, as illustrated in the following example.

\medskip
\noindent\emph{Example 1}
%\begin{example}%\label{ex:strict}
Let $n=3$, $p=2$, $\bm{\mu}=\bm{0}$, $\bm{\sigma}=\bm{1}$, $C=1$, and $\eta_j=1$, so that $\mathcal U_U^{(2)}=\{\bm{a}:|a_j|\le1,\ j=1,2,3,\ |a_1+a_2+a_3|\le\sqrt3\}$ and $\mathcal U_W^{(2)}=\{\bm{a}:\lVert\bm{a}\rVert_2\le1\}$. The vector $\bar{\bm{a}}=(1,1,-1)^\top$ belongs to $\mathcal U_U^{(2)}$ but not to $\mathcal U_W^{(2)}$, since $\lVert\bar{\bm{a}}\rVert_2=\sqrt3>1$. Its third coefficient compensates the first two in the aggregate constraint: with $a_3=0$, this constraint would be violated, since $2>\sqrt3$. For $\bm{x}=(1,1,0)^\top$, we obtain $\max_{\bm{a}\in\mathcal U_U^{(2)}}\bm{a}^\top\bm{x}=2$, attained at $\bar{\bm{a}}$, whereas $C\lVert\bm D\bm{x}\rVert_2=\sqrt2$. Hence, for $\mathcal X=\{0,1\}^3$ and $b=3/2$, $\bm{x}\in\mathcal F_W^{(2)}\setminus\mathcal F_U^{(2)}$.
%\end{example}

% Theorem~\ref{thm:dominance} applies directly to both limit theorems considered in the preceding sections. Setting $p=2$ recovers the weighted CLT formulation of Section~\ref{sec:weighted CLT}, whereas setting $p=\alpha$ recovers the weighted stable formulation of Section \ref{sec:weighted stable limit thms}. The dominance result is therefore not specific to a particular limit theorem. Rather, it follows from the norm structure induced by the scaling law of the weighted sum and the corresponding dual-norm representation of its protection term.

\section{Numerical Experiments}\label{sec:numerical comparison}
Theorem~\ref{thm:dominance} and Corollary~\ref{cor:feasible-regions} imply that, whenever $\eta_j\ge C\sigma_j$, the weighted formulation attains an optimal objective value at least as good as that of the unweighted formulation. This section quantifies the improvement and examines whether the weighted solutions attain the prescribed violation probability.

\subsection{Test problem}\label{subsec:test-problem}

We consider the $0$--$1$ knapsack problem
% \begin{equation*}\label{eq:knapsack}
%     \max_{\bm{x}\in\{0,1\}^n}\{\bm p^\top\bm x:\ \bm a^\top\bm x\le B\},
% \end{equation*}
\begin{equation*}\label{eq:knapsack}
    \max_{\bm{x}\in\{0,1\}^n}
    \Bigl\{\sum_{j=1}^n p_jx_j:\ \sum_{j=1}^n a_jx_j\le B\Bigr\},
\end{equation*}
with deterministic profits $p_j$, capacity $B$, and independent random weights $a_j$ with $\mathbb E[a_j]=\mu_j$ and $\Var(a_j)=\sigma_j^2$. Both formulations use the multiplier $C=z_{1-\varepsilon}$. Since $x_j^2=x_j$, the weighted formulation imposes
\begin{equation}\label{eq:knap weighted constraint}
    \sum_{j=1}^n\mu_jx_j+C\Bigl(\sum_{j=1}^n\sigma_j^2x_j\Bigr)^{1/2}\le B,
\end{equation}
and the unweighted formulation imposes
$\max_{\bm{a}\in\mathcal U_U}\ \sum_{j=1}^n a_jx_j\le B,$
% \begin{equation*}\label{eq:knap unweighted max}
%     \max_{\bm{a}\in\mathcal U_U}\ \sum_{j=1}^n a_jx_j\le B,
% \end{equation*} 
where $\mathcal U_U$ is the set $\mathcal U_U^{(2)}$ of Section~\ref{sec:theoretical comparison}. The parameters are chosen in favor of the unweighted formulation. First, since its aggregate bound is two-sided, asymptotic coverage $1-\varepsilon$ would require the larger multiplier $z_{1-\varepsilon/2}$ (Section~\ref{sec:weighted CLT}). Second, we set $\eta_j=C\sigma_j$, the smallest values for which Theorem~\ref{thm:dominance} applies. Since $\mathcal U_U$ is monotone in $\eta_j$, this choice is the most favorable for the unweighted formulation.

For $S=\{j:x_j=1\}$, both constraints take the form $\sum_{j\in S}\mu_j+\pi(S)\le B$, with protection terms
\begin{equation*}\label{eq:protection-terms}
    \pi_W(S)=C\Bigl(\sum_{j\in S}\sigma_j^2\Bigr)^{1/2},
    \qquad
    \pi_U(S)=C\min\Bigl\{\sum_{j\in S}\sigma_j,\
    \Bigl(\sum_{j=1}^n\sigma_j^2\Bigr)^{1/2}+\sum_{j\notin S}\sigma_j\Bigr\},
\end{equation*}
where $\pi_U(S)=\max\{\sum_{j\in S}(a_j-\mu_j):\bm{a}\in\mathcal U_U\}$. The first term inside the minimum is implied by the individual bounds on the coefficients in $S$, the second by the aggregate bound when every unselected coefficient takes its smallest value $\mu_j-C\sigma_j$, which reflects the compensation mechanism of Section~\ref{sec:theoretical comparison}. When the first term attains the minimum, the aggregate constraint is ineffective, and since $\sum_{j\in S}\sigma_j^2\le(\max_{j\in S}\sigma_j)\sum_{j\in S}\sigma_j$, we have, for $S\neq\emptyset$,
\begin{equation}\label{eq:ratio-bound}
    \frac{\pi_U(S)}{\pi_W(S)}
    =\frac{\sum_{j\in S}\sigma_j}{\bigl(\sum_{j\in S}\sigma_j^2\bigr)^{1/2}}
    \ \ge\ \frac{\bigl(\sum_{j\in S}\sigma_j^2\bigr)^{1/2}}{\max_{j\in S}\sigma_j}
    =\frac{1}{R(\bm{x})}.
\end{equation}
The reduction in protection is thus large precisely when $R(\bm{x})$ is small, which is also the regime in which the normal approximation is accurate (Proposition~\ref{prop:uniform-clt}).

\subsection{Experimental design}\label{subsec:experimental-design}

We set $\mu_j=10$ and generate $\sigma_j\sim\operatorname{Unif}\{1,\dots,5\}$ and $p_j\sim\operatorname{Unif}\{10,\dots,100\}$ independently for $j=1,\ldots,n$, so that the standard deviations range from $10\%$ to $50\%$ of the mean weights. For each $n\in\{25,50,100,250\}$, we generate $1{,}000$ instances and solve each of them for $B=10\rho n$ with $\rho\in\{0.2,0.4,0.6,0.8\}$, the fraction of items that fit if every weight equals its mean, and for $\varepsilon\in\{0.01,0.05,0.10\}$. This yields $48{,}000$ comparisons of the optimal values $z_U$ and $z_W$ of the two formulations, measured by the relative improvement $100(z_W-z_U)/z_U$. Since $\mu_j$ is constant and $\sigma_j\in\{1,\dots,5\}$, both formulations are solved exactly by a dynamic program over the number of selected items and their total standard deviation (unweighted) or variance (weighted). The experiments are implemented in Python 3.12.5 with NumPy 2.4.4 and are performed on a MacBook Pro with an Apple M3 Pro chip and 18 GB of RAM.

\subsection{Reduction in conservatism}\label{subsec:conservatism}
The weighted formulation attains a strictly larger optimal value in $47{,}886$ of the $48{,}000$ comparisons and the same value in the remaining $114$. Table~\ref{tab:density-results} reports results for $\varepsilon=0.05$.
\begin{table}[h]
\centering
\caption{Results for $\varepsilon=0.05$, averaged over $n$ (by $\rho$) and over $\rho$ (by $n$). Improvements are in \%; U and W denote selected fractions. Binding rates refer to the optimal unweighted solution, while $R(\bm x)$ is evaluated at the optimal weighted solution.}
\label{tab:density-results}
\small
\setlength{\tabcolsep}{4pt}
    \begin{tabular}{lrrrrrr}
        \toprule
        & \multicolumn{2}{c}{Improvement (\%)} & Mean & \multicolumn{2}{c}{Selected fraction} & First term\\
        \cmidrule(lr){2-3}\cmidrule(lr){5-6}
        & Mean & Median & $R(\bm x)$ & U & W & binding (\%)\\
        \midrule
        $\rho=0.2$ & 18.79 & 20.10 & 0.479 & 0.152 & 0.175 & 100.0\\
        $\rho=0.4$ & 23.29 & 24.46 & 0.333 & 0.286 & 0.357 & 100.0\\
        $\rho=0.6$ & 23.68 & 24.46 & 0.265 & 0.417 & 0.547 & 100.0\\
        $\rho=0.8$ & 17.86 & 17.75 & 0.228 & 0.582 & 0.734 & 52.1\\
        \addlinespace
        $n=25$  & 16.57 & 16.83 & 0.516 &  &  & \\
        $n=50$  & 19.87 & 20.13 & 0.366 &  &  & \\
        $n=100$ & 22.79 & 23.17 & 0.263 &  &  & \\
        $n=250$ & 24.39 & 26.09 & 0.160 &  &  & \\
        \bottomrule
    \end{tabular}
% \begin{tabular}{lrrrrrr|lrrr}
% \toprule
% & \multicolumn{2}{c}{Improvement (\%)}
% & $R(\bm x)$
% & \multicolumn{2}{c}{Selected}
% & Binding
% & & \multicolumn{2}{c}{Improvement (\%)}
% & $R(\bm x)$ \\
% \cmidrule(lr){2-3}\cmidrule(lr){5-6}
% \cmidrule(lr){9-10}
% $\rho$ & Mean & Med. & Mean & U & W & (\%)
% & $n$ & Mean & Med. & Mean \\
% \midrule
% 0.2 & 18.79 & 20.10 & 0.479 & 0.152 & 0.175 & 100.0
% & 25  & 16.57 & 16.83 & 0.516 \\
% 0.4 & 23.29 & 24.46 & 0.333 & 0.286 & 0.357 & 100.0
% & 50  & 19.87 & 20.13 & 0.366 \\
% 0.6 & 23.68 & 24.46 & 0.265 & 0.417 & 0.547 & 100.0
% & 100 & 22.79 & 23.17 & 0.263 \\
% 0.8 & 17.86 & 17.75 & 0.228 & 0.582 & 0.734 & 52.1
% & 250 & 24.39 & 26.09 & 0.160 \\
% \bottomrule
% \end{tabular}
\end{table}
For $\rho\le0.6$, the first term of $\pi_U(S)$ binds at every optimal unweighted solution, so~\eqref{eq:ratio-bound} applies, and the improvement increases with $\rho$ as the selected sets grow and $R(\bm x)$ decreases. For $\rho=0.8$, fewer unselected coefficients remain available for compensation, the second term binds in $47.9\%$ of the instances, and the improvement decreases. As $n$ increases from $25$ to $250$, the mean value of $R(\bm{x})$ decreases from $0.516$ to $0.160$, and the mean improvement increases from $16.57\%$ to $24.39\%$, again in line with~\eqref{eq:ratio-bound}. The improvement also grows with the required protection level: averaged over $n$ and $\rho$, it equals $16.43\%$, $20.90\%$, and $28.84\%$ for $\varepsilon=0.10$, $0.05$, and $0.01$, since for fixed $S$ both protection terms are linear in $C$, so that their difference grows with $C$.

{\revv \subsection{Finite-sample reliability}\label{subsec:reliability}

The weighted formulation~\eqref{eq:knap weighted constraint} is exact for independent Gaussian weights and an asymptotic approximation otherwise. To assess whether reduced conservatism compromises reliability, we set $a_j=10+\sigma_jY_j$, where the $Y_j$ are independent with zero mean and unit variance, as in Proposition~\ref{prop:uniform-clt}, and are standard normal, symmetric two-point with $\mathbb P(Y_j=1)=\mathbb P(Y_j=-1)=\tfrac12$, or shifted exponential with $Y_j=V_j-1$ and $V_j\sim\operatorname{Exp}(1)$. Since both formulations depend only on the first two moments, we evaluate the optimal solutions of Section~\ref{subsec:conservatism} for $\varepsilon=0.05$ under all three distributions, computing $\mathbb P(\bm{a}^\top\bm{x}>B)$ exactly in the Gaussian case and from $5{,}000$ Monte Carlo draws otherwise (Table~\ref{tab:reliability}).
\begin{table}[h]
    \centering
    \caption{Violation probabilities (\%) of the weighted formulation for $\varepsilon=0.05$ across $1{,}000$ instances. $R(\bm{x})$ is evaluated at the optimal weighted solution.}
    \label{tab:reliability}
    \begin{tabular}{ccrrrrrr}
        \toprule
        & & & Gaussian & \multicolumn{2}{c}{Two-point} & \multicolumn{2}{c}{Shifted exponential}\\
        \cmidrule(lr){4-4}\cmidrule(lr){5-6}\cmidrule(lr){7-8}
        $n$ & $\rho$ & Mean $R(\bm{x})$ & Mean & Mean & 95th pctl. & Mean & 95th pctl.\\
        \midrule
        25  & 0.2 & 0.739 & 3.24 & 1.50 & 6.56 & 5.21 & 6.88\\
        25  & 0.4 & 0.519 & 2.47 & 1.86 & 4.24 & 4.10 & 6.36\\
        25  & 0.6 & 0.433 & 3.75 & 3.38 & 5.02 & 5.22 & 6.66\\
        25  & 0.8 & 0.371 & 2.58 & 2.29 & 4.84 & 3.85 & 6.40\\
        \addlinespace
        50  & 0.2 & 0.527 & 3.68 & 2.77 & 4.64 & 5.31 & 6.64\\
        50  & 0.4 & 0.388 & 3.96 & 3.65 & 5.16 & 5.24 & 6.56\\
        50  & 0.6 & 0.295 & 3.87 & 3.62 & 5.02 & 5.01 & 6.26\\
        50  & 0.8 & 0.255 & 3.19 & 2.98 & 4.92 & 4.19 & 6.16\\
        \addlinespace
        100 & 0.2 & 0.406 & 4.70 & 4.38 & 5.26 & 6.00 & 6.66\\
        100 & 0.4 & 0.263 & 4.25 & 3.99 & 5.22 & 5.25 & 6.38\\
        100 & 0.6 & 0.206 & 4.35 & 4.16 & 5.22 & 5.18 & 6.24\\
        100 & 0.8 & 0.175 & 4.04 & 3.88 & 5.04 & 4.77 & 5.90\\
        \addlinespace
        250 & 0.2 & 0.243 & 4.57 & 4.31 & 5.18 & 5.49 & 6.32\\
        250 & 0.4 & 0.160 & 4.54 & 4.36 & 5.28 & 5.21 & 6.12\\
        250 & 0.6 & 0.127 & 4.08 & 3.96 & 5.00 & 4.64 & 5.68\\
        250 & 0.8 & 0.109 & 4.13 & 4.04 & 4.94 & 4.60 & 5.44\\
        \bottomrule
    \end{tabular}
\end{table}
For Gaussian weights, violation probabilities never exceed $\varepsilon$, and the means remain below $5\%$ because integrality leaves capacity slack. Since all distributions are evaluated at the same solutions, the difference between a non-Gaussian mean and the Gaussian mean is the average error of the normal approximation at the optimal weighted solutions. For every $\rho$, this error decreases monotonically as $n$ increases and $R(\bm x)$ decreases, in line with Proposition~\ref{prop:uniform-clt}. The two-point means lie below their Gaussian counterparts, whereas shifted-exponential means exceed $\varepsilon$ in 10 of the 16 settings, reaching $6.00\%$. The 95th percentiles are largest for $n=25$ and $\rho=0.2$, at $6.56\%$ and $6.88\%$, where the optimal weighted solutions select about four items and $R(\bm x)=0.74$ on average, indicating an inaccurate normal approximation, as anticipated by Corollary~\ref{cor:violation}. 
\revv{Since each violation probability has a Monte Carlo standard error of about $0.3$ percentage points, these percentiles slightly overstate the true ones.}

The unweighted formulation is markedly more conservative: averaged over $n$ and $\rho$, its violation probabilities are approximately $0.01\%$, $0.00\%$, and $0.10\%$ for the Gaussian, two-point, and shifted-exponential distributions, against $3.84\%$, $3.45\%$, and $4.95\%$ for the weighted formulation. The average improvement of $20.9\%$ at $\varepsilon=0.05$ thus stems mainly from removing unnecessary protection. When $R(\bm x)$ is large and the coefficients are strongly skewed, a slightly larger multiplier $C$ can serve as a safeguard.}

\section{Conclusion}\label{sec:conclusion}
When limit theorems are used to protect an uncertain constraint, they should be applied to the uncertain expression that appears in that constraint rather than to an unweighted aggregate. With the decision variables as weights, the weighted CLT yields the same second-order cone constraint as in the jointly Gaussian case for a much broader class of distributions, with an approximation error at any decision that is controlled by $R(\bm x)$, and weighted stable limit theorems yield power-cone constraints. Provided that the individual bounds are no tighter than $C$ standard deviations, the weighted formulation dominates its unweighted counterpart. In our knapsack experiments, it improves the optimal objective value by $20.9\%$ on average at $\varepsilon=0.05$, mainly by removing protection that is not required, with small average violation probabilities.% between $3.45\%$ and $4.95\%$.

The approach is, however, constraint-wise: joint chance constraints require a union bound, and under decision rules the weights are no longer deterministic. In these settings, the decision-independent uncertainty sets of~\cite{bandi-bertsimas} remain advantageous, since a single set protects all constraints and adaptive decisions simultaneously. %{\color{red}Natural directions for future research include the extension of the weighted approach to these settings, as well as a uniform analogue of Proposition~\ref{prop:uniform-clt} and explicit error bounds for stable domains of attraction.}
Extending the weighted approach to these cases, proving a uniform analogue of Proposition~\ref{prop:uniform-clt} and explicit error bounds for stable domains of attraction are directions for future works.

% When limit theorems are used to protect an uncertain constraint, they should be applied to the uncertain expression that appears in that constraint rather than to an unweighted aggregate. With the decision variables as weights, the weighted CLT yields the same second-order cone constraint as in the jointly Gaussian case, but for a much broader class of distributions, and the normal approximation holds uniformly over the feasible decisions when no single term dominates. Weighted stable limit theorems similarly yield power-cone constraints. Provided that the individual bounds are no tighter than $C$ standard deviations, the weighted formulation dominates its unweighted counterpart. In our knapsack experiments, it improves the optimal objective value by $20.9\%$ on average at $\varepsilon=0.05$, mainly by removing protection that is not required, while its violation probabilities remain close to the prescribed level.

% The approach is, however, constraint-wise: joint chance constraints require a union bound, and under decision rules the weights are no longer deterministic. In these settings, the decision-independent uncertainty sets of Bandi and Bertsimas~\cite{bandi-bertsimas} remain advantageous, since a single set protects all constraints and adaptive decisions simultaneously. Extending the weighted approach to these settings, and establishing a uniform analogue of Proposition~\ref{prop:uniform-clt} and explicit error bounds for stable domains of attraction, are natural directions for future research.

{\rev \section*{Statements and Declarations}
\bmhead{Funding}
The second author is funded by the NWO Rubicon grant 019.233SG.010

\bmhead{Competing interests}
The authors have no competing interests to declare.% that are relevant to the content of this article.

\bmhead{Availability of data and materials} The code to generate the instances of Sec.~\ref{subsec:experimental-design} and reproduce the results is available at \url{https://doi.org/10.5281/zenodo.23533838}.

\bmhead{Use of LLMs}
Claude (Opus 5.5) and ChatGPT (GPT-6) were used to assist with code development and the exploration of proof strategies. All outputs were reviewed and validated by the authors, who take full responsibility for the content of this manuscript.

%%===========================================================================================%%
%% If you are submitting to one of the Nature Portfolio journals, using the eJP submission   %%
%% system, please include the references within the manuscript file itself. You may do this  %%
%% by copying the reference list from your .bbl file, paste it into the main manuscript .tex %%
%% file, and delete the associated \verb+\bibliography+ commands.                            %%
%%===========================================================================================%%

\bibliography{sn-bibliography}% common bib file
\end{document}